

\baselineskip=14pt
\parskip=10pt

\magnification=\magstephalf

\def\1{{\overline{1}}}
\def\2{{\overline{2}}}
\parindent=0pt
\overfullrule=0in

\def\frac#1#2{{#1 \over #2}}
\centerline
{\bf  
The Absent-Minded Passengers Problem via Computer Algebra
}
\bigskip
\centerline
{\it Shalosh B. EKHAD and Doron ZEILBERGER}
\bigskip

{\bf Abstract.} In this {\it case study},  
we illustrate the power of {\it experimental mathematics} and {\it symbolic computation},
by {\bf discovering} interesting new facts about the so-called Absent-Minded Passengers Problem,
extending  recent work of Norbert Henze and G\"unter Last. Since we are absolutely certain
that these new facts are indeed true, and proving is not nearly as much fun as discovering,
we leave the proofs to the obtuse readers.

{\bf The Maple package.} This article is accompanied by the Maple package
{\tt AMP.txt}  that can be obtained, along with numerous input and output files, from the front of this article

{\tt http://www.math.rutgers.edu/\~{}zeilberg/mamarim/mamarimhtml/amp.html} \quad .

{\bf The initial puzzle}

The beautiful article [HL] is inspired by the following puzzle that appeared in the two delightful collections [B] and [W].

{\it An airplane with $n \geq 2$ passengers is fully booked. Passengers are boarding in chronological order, according to the
numbers on their boarding passes. The first passenger loses his boarding pass and picks one of the seats (uniformly) at random. Each
subsequent passenger takes his or her seat if available, otherwise takes one of the remaining seats (uniformly) at random.
What is the probability that the last passenger (i.e. passenger $n$) will sit in the correct seat?''}

It is not too hard to see ([B][M][HL]) that the answer is $\frac{1}{2}$. It is proved in [HL] that, more generally,
the probability  that passenger $i$ ($i \geq 2$) will sit in the correct seat is $\frac{n-i+1}{n-i+2}$.
Even more generally, they proved that when the first $k$ passengers are absent-minded, and $i>k$, that probability
equals $\frac{n-i+1}{n-i+k+1}$.

{\bf A Generatingfunctionlogy Approach to the $k=1$ case}

A quicker way to handle the original case with  only one absent-minded passenger is via {\it generating functions}
(alias {\it weight-enumerators}).

Let the {\it weight} of a sitting arrangement (a certain permutation of length $n$) resulting from
this process be the product of $w_i$ over all passengers $i$ sitting in the {\bf wrong} seat.
The initial {\it state} is when all the seats are empty. 
If, by pure luck, passenger $1$ landed in seat $1$, then the game is over, 
and the weight of that scenario is $1$ since everyone landed in the right seat.
Also the probability of that happening is $\frac{1}{n}$. Otherwise,  passenger $1$ will take seat $i$, with probability $\frac{1}{n}$,
for {\it some} $i$ between $2$ and $n$.
All the passengers, $2$ through $i-1$ will each take their rightful seat, and we now have a situation where
$i$ has to pick one of the $n-i+1$ seats in the set $\{1,i+1, \dots ,n\}$.
Let's call the initial state $S_0$ and the subsequent states $S_i$ ($2 \leq i \leq n)$).
Let $F_n(w_1, \dots, w_n)$ be the {\bf weight-enumerator} of the set  of all final sitting configurations that 
start at the initial state (our object of desire) and let $A_i$ be the weight-enumerator of  those that come from state $S_i$.

We have
$$
F_n (w_1, \dots, w_n)  = \frac{1}{n} + \frac{w_1}{n} \, \sum_{i=2}^{n} A_i \quad , \quad and
$$
$$
A_i = \frac{w_i}{n-i+1} \, \left ( 1 \, + \, \sum_{j=i+1}^{n} A_j \, \right ) \quad, \quad  2 \leq i \leq n \quad .
$$
This equation follows from the fact that passenger $i$ has $n-i+1$ equally likely choices , each of them resulting
with him sitting in the wrong seat (hence the factor $\frac{w_i}{n-i+1}$ in the front). If he chose seat $1$ then
the game is over, since all the remaining passengers seat where they are supposed to. Otherwise he sits in
seat $j$ ($i <j \leq n$), and we are in state $S_j$.

Hence, for $2 \leq i \leq n$,
$$
\frac{n-i+1}{w_i} A_i \, = \,  1 \, + \, \sum_{j=i+1}^{n} A_j \, .
$$
Replacing $i$ by $i-1$, we have
$$
\frac{n-i+2}{w_{i-1}} A_{i-1} \, = \,  1 \, + \, \sum_{j=i}^{n} A_j \, .
$$
Subtracting, we get
$$
\frac{n-i+2}{w_{i-1}} A_{i-1} - \frac{n-i+1}{w_i} A_i = A_i \quad, 
$$
implying that
$$
A_{i-1} = \frac{w_{i-1}}{w_i} \, \cdot \, \frac{w_i+n-i+1}{n-i+2} \, \cdot \, A_{i} \quad .
$$
Since $A_n=w_n$, we have, for all $0 \leq i \leq n-1$,
$$
A_{n-i}= \frac{w_{n-i}}{(i+1)!} (w_{n-i+1}+i)(w_{n-i+2}+i-1) \cdots (w_n+1) \quad .
$$
In particular (take $i=n-1$)
$$
A_1= \frac{w_1}{n!} \prod_{i=2}^{n} (w_i+n+1-i) \quad.
$$
Since $F_n(w_1, \dots, w_n)=\frac{1}{n} -\frac{w_1}{n}+A_1$, we  have:

{\bf Theorem 1}: The weight-enumerator of all sitting arrangements with one
absent-minded passenger is
$$
F_n (w_1, \dots, w_n) \, = \,
\frac{1-w_1}{n} \, + \, \frac{w_1}{n!} \prod_{i=2}^{n} (w_i+n+1-i) \quad. 
$$

It follows that for $2 \leq i \leq n$:
$$
F_n ( \, 1^{i-1} \, , \,  w_i \, , \,  1^{n-i} \,) \, = \,\frac{w_i+n+1-i}{n+2-i} \quad ,
$$
implying that the probability that passenger $i$ will sit in the {\bf right} place is $\frac{n+1-i}{n+2-i}$
(the coefficient of $w_i^0$), as proved, via a different method in [HL], Equation 1.
More generally for any subset $S$ of $\{2, \dots , n\}$, setting $w_i=1$ if $i \neq S$ and leaving $w_i$ alone when $i \in S$,
we get that the {\it marginalized} generating function equals
$$
\prod_{ i \in S} \frac{w_i+n+1-i}{n+2-i} \quad,
$$
implying that the probability that {\it all} members of $S$ will sit in the {\bf wrong} place is
$$
\prod_{ i \in S} \frac{1}{n+2-i} \quad,
$$
while the probability that they {\it all} sit in the {\bf right} place is
$$
\prod_{ i \in S} \frac{n+1-i}{n+2-i} \quad .
$$

Specializing all the $w_i$ to be $w$ we have an alternative proof of the following theorem in [HL].

{\bf Theorem 2}: The probability generating function, let's call it  $f_n(w)$, (a polynomial of degree $n$ in $w$), whose coefficient of $w^l$ is
the probability that exactly $l$ passengers sit in the wrong seat is
$$
f_n(w) \, = \, 
\frac{1-w}{n}+ \frac{w}{n!} \prod_{i=1}^{n-1} (w+i) \quad.
$$

The question of a closed-form expression for the analogous probability generating function, let's call it $f^{(k)}_n(w)$,
for the case where the first $k$ passengers are absent-minded was left open in [HL]. The next theorem
fills this gap. (Note that $f^{(1)}_n(w)=f_n(w)$).

{\bf Theorem 3}: The probability generating function $f^{(k)}_n(w)$ (a polynomial of degree $n$ in $w$), whose coefficient of $w^l$ is
the probability that exactly $l$ passengers sit in the wrong seat when the first $k$ passengers are absent-minded, is given by
$$
f^{(k)}_n(w) 
\, = \, \frac{1}{n!} \, \sum_{r=0}^{k} \, r!\,{{k} \choose {r}} w^r\,(1-w)^{k-r} \prod_{i=1}^{n-k} (rw \,+  \, i \,) \quad.
$$

We don't believe that it is possible to conjecture this theorem by merely cranking out sufficiently many special cases and
guessing a {\it pattern}. What we did was try and conjecture a generalization of Theorem 1, where one keeps track
of the actual passengers that are sitting in the wrong seat. Let $F^{(k)}_n(w_1, \dots, w_n)$ be the
multi-linear polynomial in $(w_1, \dots, w_n)$  whose coefficient of $w_{i_1} \dots w_{i_l}$ is the
probability that the passengers  in the set $\{i_1, \dots, i_l\}$ {\bf definitely} are {\it wrongly}-seated,
and the complement is {\bf definitely} seated in the {\it right} seats. (Note that $F^{(1)}_n (w_1, \dots, w_n)=F_n (w_1, \dots, w_n)$).

Using {\bf dynamical programming} (see the source code for procedure {\tt AnwkG(n,w,k)} in the Maple package {\tt AMP.txt}),
we generated lots of specific examples, that enabled 
us to discover the following generalization of Theorem 1.

{\bf Theorem 4}: Let  $e_r(w_1, \dots, w_k)$ be the coefficient of $X^r$ in $\prod_{j=1}^{k} ((1-w_j)X+w_j)$
(these are variants of the elementary symmetric functions). Then, if $n \geq k$, we have
$$
F^{(k)}_n(w_1, \dots, w_n) \, = \,
\frac{1}{n!} \sum_{r=0}^{k} r! e_{k-r}(w_1, \dots, w_k) \, \cdot \, \prod_{j=k+1}^{n} (rw_j+n+1-j) \quad .
$$

Theorem 3 follows from Theorem 4 by setting all the $w_j$'s to be $w$. Note that if we plug-in all the
$w_j$'s, {\it except} $w_i$, to be $1$, but leave $w_i$ alone, we rederive the fact, proved in [HL] another way, that
the probability of the event (if $i>k$)
`Passenger $i$  sitting in the {\bf right} seat' is $\frac{n-i+1}{n-i+k+1}$, 
Another consequence of our Theorem  4 is Theorem 3 in [HL] , that states that these events 
are independent.

By differentiating the expression for $f^{(k)}_n(w)$, given in Theorem 3, with respect to $w$, and
plugging-in $w=1$ we find (as [HL] already did) that the expectation is $k(1+\sum_{i=k+1}^{n-1} \frac{1}{i})$.
By differentiation twice, and doing some manipulatorics, one can get the expression for the variance
established in [HL]. The advantage of our Theorem 3 is that we can keep going and derive
explicit expressions for higher moments. Carsten Schneider's {\tt Sigma} package [S1][S2] should be helpful here.

{\bf The First Eight Moments of the Random Variable `Number of Passengers Sitting in the Wrong Seat' for the original case of One absent-minded passenger}

We are too lazy to find higher moments for the general case of $k$ absent-minded passengers, but we did it for
the original case of $k=1$.

Let  $X_n$ be that random variable. The expectation $E[X_n]$ , that equals $f_n'(1)$ is easily seen
(by {\it logarithmic differentiation}) to be $\sum_{i=1}^{n-1} \frac{1}{i}$, the {\it Harmonic number} $H_{n-1}$.
This is already mentioned in [HL], where they also derived an explicit expression for the variance (for arbitrary $k$),

It is convenient to introduce the notation
$$
Hn[r]:=\sum_{i=1}^{n-1} \frac{1}{i^r} \quad .
$$
Note that the upper limit is $n-1$ rather than the customary $n$. This way the formulas are much simpler.

{\bf Theorem 5}: Let $X_n$ be the random variable ``number of passengers sitting in the wrong seat" where
there is one absent-minded passenger, and $n$ passengers altogether. Then, denoting by $m_r(X_n)$ the
$r^{th}$ moment about the mean, we have (please pardon the {\it computereze})
$$
E[X_n] \,=\,  \sum_{i=1}^{n-1} \frac{1}{i} \quad .
$$
$$
Var[X_n] \,= \,(n*Hn[1]-n*Hn[2]+2*Hn[1])/n \quad .
$$
$$
m_3(X_n) \, = \,(n*Hn[1]-3*n*Hn[2]+2*n*Hn[3]-3*Hn[1]^2+6*Hn[1]-3*Hn[2])/n \quad .
$$
$$
m_4(X_n) \, = \, (3*n*Hn[1]^2-6*n*Hn[1]*Hn[2]+3*n*Hn[2]^2+4*Hn[1]^3+n*Hn[1]-7*n*Hn[2]+12*n*Hn[3] 
$$
$$
-6*n*Hn[4]-6*Hn[1]^2+14*Hn[1]-18*Hn[2]+8*Hn[3])/n \quad .
$$
$$
m_5(X_n) \, = \, (-5*Hn[1]^4+10*n*Hn[1]^2-40*n*Hn[1]*Hn[2]+20*n*Hn[1]*Hn[3]
$$
$$
+30*n*Hn[2]^2-20*n*Hn[2]*Hn[3]+10*Hn[1]^3+n*Hn[1]-15*n*Hn[2]+50*n*Hn[3]-60*n*Hn[4]
$$
$$
+24*n*Hn[5]+5*Hn[1]^2-30*Hn[1]*Hn[2]+15*Hn[2]^2+30*Hn[1]-75*Hn[2]+80*Hn[3]-30*Hn[4])/n \quad .
$$
$$
m_6(X_n) \, = \, 
(6*Hn[1]^5+15*n*Hn[1]^3-45*n*Hn[1]^2*Hn[2]+45*n*Hn[1]*Hn[2]^2-15*n*Hn[2]^3-15*Hn[1]^4
$$
$$
+25*n*Hn[1]^2-180*n*Hn[1]*Hn[2]+220*n*Hn[1]*Hn[3]-90*n*Hn[1]*Hn[4]+195*n*Hn[2]^2
$$
$$
-300*n*Hn[2]*Hn[3]+90*n*Hn[2]*Hn[4]+40*n*Hn[3]^2+20*Hn[1]^3+n*Hn[1]-31*n*Hn[2]
$$
$$
+180*n*Hn[3]-390*n*Hn[4]+360*n*Hn[5]-120*n*Hn[6]+90*Hn[1]^2-330*Hn[1]*Hn[2]
$$
$$
+120*Hn[1]*Hn[3]+225*Hn[2]^2
$$
$$
-120*Hn[2]*Hn[3]+62*Hn[1]-270*Hn[2]+520*Hn[3]-450*Hn[4]+144*Hn[5])/n \quad .
$$

For $m_7(X_n)$ and $m_8(X_n)$ see the web-page

{\tt https://sites.math.rutgers.edu/\~{}zeilberg/tokhniot/oAMP2a.txt} \quad .

{\bf Sketch of the Proof}: The above theorem was discovered by {\bf pure guessing}, using
an {\it ansatz} with {\it undetermined coefficients} featuring the quantities
$Hn[r]$, that are the partial sums of $\zeta(r)$ for $r \geq 2$. They are
all  solutions of complicated recurrences and hence can be rigorously proved using
Carsten Schneider's amazing {\tt Sigma} package [S1][S2]. Since we are sure that they are true,
we did not bother to actually do it.

The web-page also has asymptotic expansions for these quantities, confirming, via elementary means,
that $X_n$ is {\bf asymptotically normal}, up to the $8^{th}$ moment, and one can easily go far beyond.
This fact was proved using `advanced' probability in [HL].

{\bf Recurrences for $f^{(k)}_n(w)$ for $k=1,2,3,4$}

It is useful to have recurrences for these quantities.

{\bf Theorem 6}: The probability generating function for the random variable,  
`the number of passengers sitting in the wrong seat where the first $k$ passengers are absent-minded',
$f^{(k)}_n(w)$ satisfy the following linear recurrences.

For $ f^{(1)}_n(w)$ (alias $f_n(w)$) we have
$$
{\frac {n \left( n+w \right) f^{(1)}_n(w) } { \left( 2+n \right)  \left( 1+n \right) }}-{\frac { \left( 2\,n+w+1 \right) 
 f^{(1)}_{n+1}(w) }{2+n}}+ f^{(1)}_{n+2}(w) =0
$$
For $f^{(2)}_n(w)$ we have
$$
-{\frac {n \left( n+2\,w \right)  \left( n+w \right) f^{(2)}_n(w)}{ \left( n+4 \right)  \left( n+3 \right)  \left( 2+n \right) }}
+{\frac { \left( 3\,{n}^{2}+6\,nw+2\,{w}^{2}+3\,n+3\,w+1 \right)  f^{(2)}_{n+1}(w) }{ \left( n+4 \right)  \left( n+3 \right) }} \quad ,
$$
$$
-3\,{\frac { \left( n+w+1
 \right)  f^{(2)}_{n+2}(w) }{n+4}}+ f^{(2)}_{n+3}(w) =0 \quad .
$$
For $f^{(3)}_n(w)$ we have
$$
{\frac {n \left( n+2\,w \right)  \left( n+3\,w \right)  \left( n+w \right) f^{(3)}_n(w) }{ \left( n+5 \right)  \left( n+4 \right)  \left( n+3
 \right)  \left( n+6 \right) }}
-{\frac { \left( 3\,w+1+2\,n \right)  \left( 2\,{n}^{2}+6\,nw+2\,{w}^{2}+2\,n+3\,w+1 \right) f^{(3)}_{n+1}(w) }{
 \left( n+5 \right)  \left( n+4 \right)  \left( n+6 \right) }}
$$
$$
+{\frac { \left( 6\,{n}^{2}+18\,nw+11\,{w}^{2}+12\,n+18\,w+7 \right)  f^{(3)}_{n+2}(w) }
{ \left( n+6 \right)  \left( n+5 \right) }}-2\,{\frac { \left( 2\,n+3\,w+3 \right)  f^{(3)}_{n+3}(w) }{n+6}}+ f^{(3)}_{n+4}(w) =0
\quad .
$$
For a recurrence for $f^{(4)}_{n}(w)$ see the web-page

{\tt https://sites.math.rutgers.edu/\~{}zeilberg/tokhniot/oAMP1.txt} \quad .

One can easily go further. In general $f^{(k)}_n(w)$ satisfies a linear recurrence equation of order $k+1$.

{\bf References}

[B] Bela Bolob\'as, {\it ``The Art of Mathematics: Coffee Time in Memphis''}, Cambridge University Press, 2006.

[HL] Norbert Henze and G\"unter Last, {\it Absent-Minded Passengers}, American Mathematical Monthly {\bf 126}(10) (Dec. 2019),
867-875.

[S1] Carsten Schneider, {\it The Summation package Sigma}, A Mathematica package available from \hfill\break
{\tt https://www3.risc.jku.at/research/combinat/software/Sigma/index.php} 

[S2] Carsten Schneider, {\it  Symbolic Summation Assists Combinatorics}, Sem.Lothar.Combin. {\bf 56}(2007),Article B56b (36 pages). \hfill\break
{\tt https://www3.risc.jku.at/research/combinat/software/Sigma/pub/SLC06.pdf}

[W] Peter Winkler, {\it ``Mathematical Puzzles: A Connoisseur's Collection''}. A.K. Peters, 2004.

\bigskip
\bigskip
\hrule
\bigskip
Shalosh B. Ekhad, c/o D. Zeilberger, Department of Mathematics, Rutgers University (New Brunswick), Hill Center-Busch Campus, 110 Frelinghuysen
Rd., Piscataway, NJ 08854-8019, USA. \hfill\break
Email: {\tt ShaloshBEkhad at gmail dot com}   \quad .
\bigskip
Doron Zeilberger, Department of Mathematics, Rutgers University (New Brunswick), Hill Center-Busch Campus, 110 Frelinghuysen
Rd., Piscataway, NJ 08854-8019, USA. \hfill\break
Email: {\tt DoronZeil at gmail  dot com}   \quad .
\bigskip
\hrule
\bigskip
{\bf Exclusively published in the Personal Journal of Shalosh B.  Ekhad and Doron Zeilberger and arxiv.org \quad .}
\bigskip
\hrule
\bigskip
Written: Jan. 19, 2020.
\end